\documentclass{article}%
\pdfoutput=1
\usepackage{amsmath}
\usepackage{amsfonts}
\usepackage{amssymb}
\usepackage{graphicx}%
\providecommand{\U}[1]{\protect\rule{.1in}{.1in}}
\providecommand{\U}[1]{\protect\rule{.1in}{.1in}}
\newtheorem{theorem}{Theorem}
\newtheorem{lemma}[theorem]{Lemma}

\begin{document}

\title{A Simple Proof of Thue's Theorem \\on Circle Packing}
\author{Hai-Chau Chang\thanks{Retired from Department of Mathematics, National Taiwan
University} and Lih-Chung Wang\thanks{Corresponding author, Department of
Applied Mathematics, National Donghwa University, Shoufeng, Hualien 974
Taiwan, R.O.C. Email: lcwang@mail.ndhu.edu.tw} \thanks{Paper partially
supported by National Science Council (NSC-98-2115-M-259-003)} }
\date{}
\maketitle

%%\begin{abstract}
%%A simple proof of Thue theorem on Circle Packing is Given. The proof is only
%%based on density analysis of Delaunay triangulation for the set of points that
%%are centers of circles in a saturated circle configuration.
%%\end{abstract}

Thue's theorem states that the regular hexagonal packing is the densest circle
packing in the plane. The density of this circle configuration is%
\[
\frac{\pi}{\sqrt{12}}\approx0.906\,90.
\]

In geometry, circle packing refers to the study of the arrangement of unit
circles on the plane such that no overlapping occurs, which is the
2-dimensional analog of Kepler's sphere packing problem proposed in 1611. A
\textit{circle configuration} which refers to the centers of circles is a set
of points such that the distance between any two points in the set is greater
than or equal to $2$. Imagine filling a large container with small unit
circles inside. The density of the arrangement is the proportion of the area
of the container that is taken up by the circles. In order to maximize the
number of circles in the container, you need to find an arrangement with the
highest possible density, so that the circles are packed together as closely
as possible. Hence, the density of a circle configuration is the asymptotic
limit on density with the container getting bigger and bigger. In 1773,
Lagrange proved that the minimal density is $\pi/\sqrt{12}$ by assuming that
the circle configurations are lattices. In 1831, Gauss proved that the minimal
density of sphere packing is $\pi/\sqrt{18}$ by assuming that the sphere
configurations are lattices. Without the lattice assumption, the first proof
of circle packing problem was made by Axel Thue. However, it is generally
believed that Thue's original proof was incomplete and that the first complete
and flawless proof of this fact was produced by L. F. Toth (1940). Later,
different proofs were proposed by Segre and Mahler, Davenport, and Hsiang.

A circle configuration is called \textit{saturated }if it is not a proper
subset of another circle configuration. Given a circle configuration
$\mathcal{C}$,\ any saturated circle configuration containing $\mathcal{C}%
$\ \ is called a \textit{saturation} of $\mathcal{C}$. Since the density of a
circle configuration $\mathcal{C}$ is always less than or equal to the density
of any saturation of $\mathcal{C}$. Hence, we only need to consider the
saturated circle configurations instead of all circle configurations.%
%TCIMACRO{\FRAME{dtbpFU}{4.9502in}{5.0652in}{0pt}{\Qcb{A Delaunay triangulation
%in the plane with circumcircles shown. From
%http://en.wikipedia.org/wiki/Delaunay\_triangulation}}{}%
%{delaunay_circumcircles.png}{\special{ language "Scientific Word";
%type "GRAPHIC";  maintain-aspect-ratio TRUE;  display "USEDEF";
%valid_file "F";  width 4.9502in;  height 5.0652in;  depth 0pt;
%original-width 4.8957in;  original-height 5.0107in;  cropleft "0";
%filename 'Delaunay_circumcircles.png';file-properties "XNPEU";}}}%
%BeginExpansion
\begin{center}
\includegraphics[
natheight=5.010700in,
natwidth=4.895700in,
height=5.0652in,
width=4.9502in
]%
{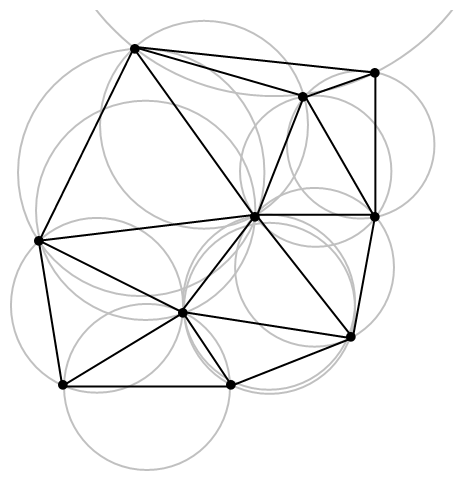}%
\\
A Delaunay triangulation in the plane with circumcircles shown. From
http://en.wikipedia.org/wiki/Delaunay\_triangulation
\end{center}
%EndExpansion

In computational geometry, a point set triangulation, i.e., a triangulation
$T(\mathcal{C})$\ of a discrete set of points $\mathcal{C}$\ on the plane is a
subdivision of the convex hull of the points into triangles such that any two
triangles intersect in a common edge or not at all and the set of points that
are vertices of the triangles coincides with $\mathcal{C}$\ . The Delaunay
triangulation for a set $\mathcal{C}$\ of points in the plane is a
triangulation $DT(\mathcal{C})$ such that no point in the set $\mathcal{C}%
$\ is inside the circumcircle of any triangle in $DT(\mathcal{C})$. Delaunay
invented such triangulations in 1934. The uniqueness and existence of Delaunay
triangulations are both uncertain. For example, there is no Delaunay
triangulation for a set of points on a straight line. The Delaunay
triangulations for four points on a circle are not unique; it is obvious that
there are two possible triangulations for a cocircular quadrilateral splitting
into two triangles. However, there always exists a Delaunay triangulation for
a saturated circle configuration. To find the Delaunay triangulation of a set
of points in the plane can be converted to find the convex hull of a set of
points in $3$-dimensional Euclidean space, by giving every point $p$ in a
saturated circle configuration an extra coordinate equal to $|p|^{2}$, taking
the convex hull, and mapping back to the Euclidean plane by forgetting the
last coordinate. A facet of the convex hull not being a triangle implies that
at least $4$ of the original points lay on the same circle, which makes the
triangulation not unique.

\begin{lemma}
Let $\theta$ be the largest internal angle of a triangle $\Delta ABC$ in a
Delaunay triangulation for a saturated circle configuration $\mathcal{C}$.
Then%
\[
\frac{\pi}{3}\leq\theta<\frac{2\pi}{3}.
\]

\end{lemma}

\noindent\textbf{Proof:} The largest internal angle of a triangle is always
bigger than or equal to $\frac{\pi}{3}$. The equality only holds for regular triangles.

Suppose that $\theta\geq\frac{2\pi}{3}$. Let say $A$ to be the smallest
internal angle. We have $\sin A\leq\frac{1}{2}$ and $\overline{BC}\geq2$.
Denote the circumradius of $\Delta ABC$ by $R$. By the sine law, we have
\[
2R=\frac{\overline{BC}}{\sin A}\geq\frac{2}{\sin A}\geq4.
\]
Then the circumcenter of $\Delta ABC$ can be added to the circle configuration
$\mathcal{C}$ which contradicts the saturated-ness of the circle configuration
$\mathcal{C}$. Therefore, we obtain%
\[
\theta<\frac{2\pi}{3}. \,\,\,\,\mbox{ Q.E.D.}
\]

The density of a triangle $\Delta ABC$ in a Delaunay triangulation for a
saturated circle configuration $\mathcal{C}$ is equal to
\[
\frac{\frac{1}{2}A+\frac{1}{2}B+\frac{1}{2}C}{\text{ the area of }\Delta
ABC}=\frac{\pi/2}{\text{ the area of }\Delta ABC}.
\]

\bigskip

\begin{lemma}
The density of a triangle $\Delta ABC$ in a Delaunay triangulation for a
saturated circle configuration $\mathcal{C}$ is less than or equal to
$\pi/\sqrt{12}$. The equality holds only for the regular triangle with
side-length $2$.
\end{lemma}

\noindent\textbf{Proof:} Let say that $B$ is the largest internal angle of
$\Delta ABC$. Then, by the above lemma,
\[
\text{the area of }\Delta ABC=\frac{1}{2}\overline{AB}\cdot\overline{BC}%
\cdot\sin B\geq\frac{1}{2}\cdot2\cdot2\cdot\frac{\sqrt{3}}{2}=\sqrt{3}.
\]
Therefore, we have
\[
\text{the density of }\Delta ABC=\text{ }\frac{\pi/2}{\text{ the area of
}\Delta ABC}\leq\frac{\pi}{\sqrt{12}}.
\]
It is obvious from the computation that the equality holds only when $\Delta
ABC$ is a regular triangle and side-length of $\Delta ABC$ is $2$. Q.E.D.

The density of the union of any finite Delaunay triangles in a saturated
circle configuration is a weighted average of the densities of the Delaunay
triangles. i.e.
\[
\text{the density }=\frac{\sum_{\Delta_{i}\text{:Delaunay triangle}}(\text{the
area of }\Delta_{i})\times(\text{the density of }\Delta_{i})}{\sum_{\Delta
_{i}\text{:Delaunay triangle}}\text{the area of }\Delta_{i}}.
\]
Since we have shown that the density of a Delaunay triangle is less than or
equal to $\pi/\sqrt{12}$, the density of the union of any finite Delaunay
triangles in a saturated circle configuration is also less than or equal to
$\pi/\sqrt{12}$. Therefore, we obtain a simple proof of Thue theorem.

\begin{theorem}
[Axel Thue]The hexagonal lattice is the densest of all possible circle packings.
\end{theorem}

\end{document}